\documentclass[preprint,3p,sort]{elsarticle}

\usepackage{listings}
\usepackage{amsfonts}
\usepackage{amsthm}
\usepackage{graphicx}
\usepackage{epstopdf}
\usepackage{algorithmic}
\usepackage{amsmath}
\usepackage{amssymb}
\usepackage{hyperref}
\usepackage[T1]{fontenc}
\usepackage{soul,color}
\newtheorem{theorem}{Theorem}[section]
\newtheorem{proposition}[theorem]{Proposition}

\theoremstyle{remark}
\newtheorem{remark}[theorem]{Remark}
\theoremstyle{example}

\numberwithin{equation}{section}

\begin{document}

\newcommand{\Surf}{{\partial \Omega}}
\newcommand{\Vol}{{\Omega}}
\newcommand{\dVol}{{\Gamma}}
\newcommand\txtred[1]{{\color{red}#1}}
\newcommand\txtblue[1]{{\color{blue}#1}}
\newcommand{\xj}{{x_j}}
\newcommand{\ui}{{u_i}}
\newcommand{\uj}{{u_j}}
\newcommand{\pd}[2]{\frac{\partial #1}{\partial #2}}
\newcommand{\taupr}{\tau^\prime}
\newcommand{\Rot}{\mathbb{R}}

\newcommand{\half}{\frac{1}{2}}
\newcommand{\BigO}[1]{\ensuremath{\operatorname{O}\bigl(#1\bigr)}}
\newcommand{\Order}{\mathcal{O}}
\newcommand{\U}{{\mathbf{U}}}
\newcommand{\Ue}{{\mathbf{V}}}
\newcommand{\Ux}{{\mathbf{U}_x}}
\newcommand{\Vv}{{\mathbf{V}}}
\newcommand{\vv}{{\overline{v}}}
\newcommand{\G}{{\mathbf{G}}}
\newcommand{\At}{\tilde{A}}
\newcommand{\Pp}{\mathcal{P}}
\newcommand{\Ppo}{{\mathcal{\overline{P}}}}
\newcommand{\Q}{\mathcal{Q}}
\newcommand{\Qx}{\mathcal{Q}_{x}}
\newcommand{\Qxe}{\mathcal{Q}_x^e}
\newcommand{\Qxx}{\mathcal{Q}_{xx}}
\newcommand{\Qxi}{\mathcal{Q}_{x^i}}
\newcommand{\B}{\mathcal{B}}
\newcommand{\Lagr}{\mathcal{L}}
\newcommand{\lagr}{\mathcal{l}}
\newcommand{\D}{\mathcal{D}}
\newcommand{\Dx}{\mathcal{D}_{x}}
\newcommand{\Dxi}{\mathcal{D}_{x^i}}
\newcommand{\Dxa}{\mathcal{D}_{x^1}}
\newcommand{\Dxb}{\mathcal{D}_{x^2}}
\newcommand{\Dxc}{\mathcal{D}_{x^3}}
\newcommand{\Bx}{\B\Dx}
\newcommand{\Lx}{\mathbb{L}_x}
\newcommand{\Di}{\mathcal{D}_{AD}}
\newcommand{\Dii}{\mathcal{D}_{AD_i}}
\newcommand{\De}{\mathcal{D}^{e}}
\newcommand{\Dxx}{\mathcal{D}_{xx}}
\newcommand{\Ee}{\mathcal{E}}
\newcommand{\R}{\mathcal{R}}
\newcommand{\nn}{\nonumber}
\newcommand{\Nn}{\ell}
\newcommand{\Nns}{\ell^\star}
\newcommand{\Mm}{\mathcal{M}}
\newcommand\norm[1]{\left\lVert#1\right\rVert}
\newcommand\bint[3]{\left. #1 \right|_{#2}^{#3}}
\newcommand\sg[1]{\sigma^#1}
\newcommand\eps{\varepsilon}
\newcommand\epst{\overline{\varepsilon}}
\newcommand{\V}{\overline{u}}
\newcommand{\F}{\mathcal{F}}
\newcommand{\Aa}{\mathcal{A}}
\newcommand\bb[1]{\textcolor{blue}{#1}}
\newcommand\rred[1]{\textcolor{red}{#1}}
\newcommand\ggreen[1]{\textcolor{green}{#1}}
\newcommand\bblue[1]{\textcolor{blue}{#1}}
\newcommand{\I}{\mathcal{\tilde{I}}}
\newcommand{\A}{\mathcal{\tilde{A}}}
\newcommand{\Iha}{\mathcal{\hat{I}}_{4t2}}
\newcommand{\Ihb}{\mathcal{\hat{I}}_{4t4}}
\newcommand{\al}{\alpha}
\newcommand\bft[1]{{\bf{#1}}}
\newcommand{\drho}{\Delta\rho}
\newcommand{\BT}{\mathbb{BT}}
\newcommand{\BTI}{\mathbb{BTI}}
\newcommand{\SAT}{\mathbb{SAT}}
\newcommand{\BC}{\mathbb{BC}}
\newcommand{\Ut}{\tilde{U}}
\newcommand{\C}{\mathbb{C}}

\makeatletter
\newcommand{\doublewidetilde}[1]{{%
  \mathpalette\double@widetilde{#1}%
}}
\newcommand{\double@widetilde}[2]{%
  \sbox\z@{$\m@th#1\widetilde{#2}$}%
  \ht\z@=.9\ht\z@
  \widetilde{\box\z@}%
}
\makeatother


\begin{frontmatter}

%

\title{An Energy Stable Incompressible  Multi-Phase Flow Formulation}

\author[sweden,southafrica]{Jan Nordstr\"{o}m}
\corref{firstcorrespondingauthor}
\cortext[firstcorrespondingauthor]{Corresponding author}
\ead{jan.nordstrom@liu.se}
\author[capetown]{Arnaud.G. Malan}
\address[sweden]{Department of Mathematics,  Applied Mathematics, Link\"{o}ping University, SE-581 83 Link\"{o}ping, Sweden}
\address[southafrica]{Department of Mathematics   $\&$ Applied Mathematics, University of Johannesburg, Auckland Park 2006, South Africa}
\address[capetown]{InCFD Research Group,  Department of Mechanical Engineering, University of Cape Town, Cape Town 7700, South Africa}

\begin{abstract}
We show that a reformulation of the governing equations for incompressible multi-phase flow in the volume of fluid setting leads to a well defined energy rate. New nonlinear inflow-outflow and solid wall boundary conditions bound the energy rate and lead to an energy estimate in terms of only external data. The new formulation combine perfectly with summation-by-parts operators and leads to provable energy stability. 
\end{abstract}

\begin{keyword}
multi-phase flow \sep  volume of fluid
\sep boundary conditions \sep energy stability \sep summation-by-parts

\end{keyword}


\end{frontmatter}


\section{Introduction}

Initial boundary value problems (IBVPs) for nonlinear flow problems including boundary conditions are notoriously difficult to bound.
We have previously  \cite{nordstrom2022linear-nonlinear,Nordstrom2022_Skew_Euler,NORDSTROM2024_BC,nordstrom2024skewsymmetric_jcp}, reformulated the shallow water equations and the incompressible and compressible Euler and Navier-Stokes equations 
such that energy estimates
were obtained. 
We also discretized the new formulations and arrived at provably stable nonlinear schemes \cite{nordstrom_roadmap}.  


In this note we provide a theoretical background for energy stability of the IBVP for incompressible multi-phase liquid-gas flows in the volume-of-fluid (VOF) formulation  \cite{Hirt1981}. Specific modeling techniques for sharpening and diffusing the interface are for now, left to others \cite{WOS:000284670300010,Mani}.
The VOF formulation is applicable to complex interface motions, it is mass conservative and tracks the interface 
by advecting the volume fraction of the target phase.  Combined with a single liquid-gas velocity, a "one-fluid" formulation \cite{Oxtoby2015,Malan2021} results. 

In \cite{Liu2007,Shen2015279,YANG2019229,PAN2023107329}, the original equations are similar, but the dependent variables and bounds differ and most importantly: boundary conditions are essentially ignored. Here we reformulate the one-fluid VOF equations into a new set of skew-symmetric equations and derive new boundary conditions that lead to an energy bound. By discretizing using summation-by-parts (SBP) operators, energy stability follows.

 \section{The reformulation}
\label{sec:CF}
We consider an  incompressible viscous liquid ($l$) and gas ($g$) mixture in two dimensions (2D)
 (with trivial extension to 3D). Using  Einstein's summation convention, the classical one-fluid VOF formulation  \cite{Hirt1981}  reads
\begin{align}\label{eq:ge0_a}
\partial_t \alpha + u_j \partial_j  \alpha & = 0	 \nonumber \\
\partial_t u_i + u_j  \partial_j  u_i & = \partial_j \tau_{ij}/\rho,  \quad i=1,2\\
\partial_j u_{j} & = 0.  \nonumber
\end{align}
In (\ref{eq:ge0_a}), $\alpha$ is the volume fraction of the liquid, $u_j$ is velocity in direction $x_j$, $\tau_{ij}  = \tau_{ij}^* - \delta_{ij} p$ is  the stress tensor,   $\tau_{ij}^*  = \left[ \mu \left(u_{i,j} + u_{j,i} \right) \right] $ is the viscous stress tensor and $p $ is pressure. Furthermore,  $\rho  = \alpha \rho_l+(1-\alpha)\rho_g$ and $\mu = \alpha \mu_l+(1-\alpha)\mu_g$  are the volume-averaged density and viscosity, respectively. The derivatives are denoted $\partial_t \psi=\psi_{,t}=\partial \psi/ \partial t$ and $\partial_j  \psi=\psi_{,j}=\partial \psi / \partial x_j$, We neglected the external gravity forces which have no impact on stability. To get at an energy bound, the formulation (\ref{eq:ge0_a}) will be modified in two steps. 

In the first step we replace the volume fraction by the density as dependent variable 
 and move the divergence relation to the righthand side leading to ($S_{ij}= \tau_{ij}$ and $S_{3j}= -u_{j}$) the equivalent formulation
\begin{align}\label{eq:ge1}
\partial_t \rho + u_j \partial_j  \rho & = 0	 \nonumber \\
\partial_t u_i + u_j  \partial_j u_i & = \partial_j S_{ij}/\rho,  \quad i=1,2\\
0 & = \partial_j S_{3j}.  \nonumber
\end{align}

In the second step we aim for a scaling of the viscous terms and introduce the new variables $[ \phi_0, \phi_j, \phi_3]^T=[\sqrt{\rho},  \sqrt{\rho} u_j, p]$ (see  \cite{nordstrom2022linear-nonlinear,Nordstrom2022_Skew_Euler,NORDSTROM2024_BC,nordstrom2024skewsymmetric_jcp}, for similar but not identical choices) into (\ref{eq:ge1}) to yield the new equation set
\begin{align}\label{eq:ge2}
\partial_t \phi_0 + u_j \partial_j  \phi_0 & = 0	 \nonumber \\
\partial_t \phi_i + u_j  \partial_j \phi_i & = \partial_j S_{ij}/\phi_0,  \quad i=1,2\\
0 & = \partial_j S_{3j}.  \nonumber
\end{align}

Introducing $\Phi=[\phi_0,\phi_1,\phi_2,\phi_3]^T$,  $\tilde I =diag[1,1,1,0]$, $\tilde A_j =diag[u_j, u_j, u_j,0]$, $\tilde B =diag[0, 1/\phi_0,1/\phi_0,1]$, ${S_j}=(0,S_{1j} ,S_{2j},S_{3j})^T$ and noting that $\tilde A_{j,j}$ vanish, we  cast  (\ref{eq:ge2}) in the final  matrix-vector form
\begin{equation} \label{eq:ge4}
\tilde I  \partial_t  \Phi  + \frac{1}{2} \Big[  \partial_j  (\tilde A_j \Phi) + \tilde A_j   \partial_j \Phi  \Big] =  \tilde B \partial_j S_j.
\end{equation}
\begin{remark} \label{Rem1}
The first step leading to (\ref{eq:ge1}) produces the skew-symmetric lefthand side of (\ref{eq:ge4}) suitable for Green's theorem. The second step leading to (\ref{eq:ge2}) produces a scaling of  the righthand side in (\ref{eq:ge4}) again opening up for Green's theorem. Reformulating (\ref{eq:ge0_a}) into (\ref{eq:ge4}) simplifies  the upcoming energy analysis.
\end{remark}

\section{Energy analysis}
Focusing first on the skew-symmetric lefthand side of  (\ref{eq:ge4}), we multiply with $\Phi^T$ and integrate to get
\begin{equation} \label{eq:ener1}
\int_\Omega \Phi^T \tilde I  \partial_t  \Phi  d\Omega + \frac{1}{2}  \int_\Omega \Big[ \Phi^T  \partial_j  (\tilde A_j \Phi) + \Phi^T \tilde A_j   \partial_j \Phi  \Big] d\Omega = \int_\Omega \Phi^T  \tilde B \partial_j S_j d\Omega.
\end{equation}
Letting $\int_\Omega  \Phi^T \tilde I \Phi d\Omega = ||\Phi||_{\tilde I}^2$ and $ \Phi^T  \partial_j  (\tilde A_j \Phi) =\partial_j (\Phi^T \tilde A_j \Phi) - (\partial_j \Phi)^T \tilde A_j \Phi$ followed by Green's theorem yields
\begin{equation} \label{eq:ener2}
 \frac{d}{dt}  ||\Phi||_{\tilde I}^2 +  \oint\limits_{\partial\Omega}\Phi^T \left( n_j \tilde A_j \right) \Phi ds = 2  \int_\Omega \Phi^T  \tilde B \partial_j S_j  d\Omega.
\end{equation}
In (\ref{eq:ener2}), $\partial\Omega$, $ds$ and $n_j$ are respectively the boundary, its surface element and outward pointing unit normal.  

Focusing secondly on the rescaled righthand-side of (\ref{eq:ge4}), we see that the righthand side of (\ref{eq:ener2}) implies
\begin{align} \label{eq:enerS}
  \int_\Omega  \Phi^T  \tilde B \partial_j S_j  d\Omega & =  \int_\Omega \Big[ \phi_i  \partial_j S_{ij}/\phi_0 + \phi_3  \partial_j S_{3j} \Big] d\Omega = \int_\Omega \Big[ u_i \tau_{ij,j} - p u_{i,i}\Big] d\Omega =
  \int_\Omega  \Big[ (u_i \tau_{ij})_{,j} -  u_{i,j}  \tau_{ij}- p u_{i,i} \Big] d\Omega \nonumber \\ & =  \int_\Omega  \Big[ (u_i \tau_{ij})_{,j} -  u_{i,j} \left( \tau^*_{ij}- \delta_{i,j} p \right)- p u_{i,i} \Big] d\Omega  =  \oint\limits_{\partial\Omega} u_i \tau_{ij} n_j ds -  \int_\Omega u_{i,j} \tau^*_{ij} d\Omega,
\end{align}
where we again used Green's theorem.  
Inserting (\ref{eq:enerS}) into (\ref{eq:ener2}) and rearranging leads to the energy rate
\begin{equation} \label{eq:ener3}
 \frac{d}{dt}  ||\Phi||_{\tilde I}^2  + 2  \int_\Omega  u_{i,j} \tau^*_{ij} d\Omega +  \oint\limits_{\partial\Omega}  \BT ds = 0,
\end{equation}
where the $2^{nd}$ term ($u_{i,j} \tau^*_{ij} = 2 u_{i,i}^2 + (u_{1,2} + u_{2,1})^2 \geq 0$) provides dissipation and the  boundary term $\BT$ is
\begin{align} \label{eq:enerBT2}
\BT = \Phi_1^T A \Phi_1 - 2 \left[ u_n (\tau_n -p) + u_t \tau_{t} \right].
\end{align}
In (\ref{eq:enerBT2}), $A=diag[u_n,u_n,u_n]$ and $\Phi_1^T=[\phi_0, \phi_n, \phi_t]=[ \phi_0, \phi_0 u_n, \phi_0 u_t]$.
The subscripts $n,t$ denote components normal and tangential to $\partial\Omega$. 
For the stress term in (\ref{eq:enerBT2}) we use  ${\bf{\tau}}^T = [\tau_1,\tau_2]=[ \tau^*_{1j}n_j, \tau^*_{2j}n_j]$ and the rotation matrix $N$ (given in (\ref{Def_R1_N} below) to obtain 
${\bf{u}}^T   {\bf{\tau}} = {\bf{u}}^T \left( N^T N \right)  {\bf{\tau}}  = (N {\bf{u}})^T \left( N  {\bf{\tau}}  \right)= [u_n,u_t] \cdot [\tau_n ,\tau_t]^T$.

\subsection{The boundary conditions}
The unrestrained boundary term  in (\ref{eq:enerBT2}) can be written in vector-matrix-vector form as
\begin{equation}\label{eq:enerBT3}
\BT = \begin{pmatrix}  \Phi_1 \\ \Phi_2 \end{pmatrix}^T
\begin{pmatrix} A & B \\  B^T & 0 \end{pmatrix} 
\begin{pmatrix}  \Phi_1 \\ \Phi_2 \end{pmatrix}, \quad \text{where} \quad 
 \Phi_2 = \begin{pmatrix} p \\ \tau_n \\ \tau_t \end{pmatrix} \quad \text{and} \quad 
B=\frac{1}{\phi_0} \begin{pmatrix}  0& 0 & 0 \\ 1 &-1 & 0 \\0 & 0 & -1
\end{pmatrix}
\end{equation} 
while $\Phi_1$ and $A$ are given above.
Next, we introduce the non-singular block rotation matrix $\Rot$ as
\begin{equation}\label{eq:rotCont}
\Rot \Rot^{-1} = \begin{pmatrix} I & R \\ 0 & I \end{pmatrix}
\begin{pmatrix} I & -R \\ 0 & I \end{pmatrix} = \begin{pmatrix} I & 0 \\ 0 & I \end{pmatrix},
\end{equation}
where the blocks are $3 \times 3$ matrices. The rotation matrix (\ref{eq:rotCont}) introduced into (\ref{eq:enerBT3}) yields
\begin{align}\label{eq:enerBT4}
\BT & = \begin{pmatrix} \Phi_1 \\ \Phi_2 \end{pmatrix}^T \Rot \Rot^{-1} 
\begin{pmatrix} A& B \\  B^T & 0 \end{pmatrix} \Rot \Rot^{-1} \begin{pmatrix}  \Phi_1 \\ \Phi_2 \end{pmatrix} \nn \\
& = \left( \Rot^{-1} \begin{pmatrix} \Phi_1 \\ \Phi_2 \end{pmatrix} \right)^T \begin{pmatrix} A & A R + B \\  R^T A + B^T & R^T (A R + B) +B^T R \end{pmatrix}  \left(  \Rot^{-1} \begin{pmatrix}  \Phi_1 \\ \Phi_2 \end{pmatrix} \right).
\end{align}
With matrix $A$ being non-singular, we cancel the off-diagonal matrices in (\ref{eq:enerBT4}) with $R=-A^{-1}B$ to get
\begin{align}\label{eq:enerBT5}
\BT &= \begin{pmatrix} A^{-1} [ A \Phi_1 + B \Phi_2] \\ \Phi_2 \end{pmatrix}^T \begin{pmatrix} A & 0 \\  0 & -B^T A^{-1} B \end{pmatrix} \begin{pmatrix} A^{-1} [ A \Phi_1 + B \Phi_2] \\ \Phi_2 \end{pmatrix} \nn \\
& =  \begin{pmatrix} A \Phi_1 + B \Phi_2 \\ B \Phi_2 \end{pmatrix}^T \begin{pmatrix} I/u_n & 0 \\  0 & -I/u_n \end{pmatrix} \begin{pmatrix} A \Phi_1 + B \Phi_2 \\ B \Phi_2 \end{pmatrix} =  \begin{pmatrix} W_1 \\ W_2 \end{pmatrix}^T \begin{pmatrix} \Lambda_1 & 0 \\  0 & \Lambda_2 \end{pmatrix} 
  \begin{pmatrix} W_1 \\ W_2 \end{pmatrix},  \  \text{where}
\end{align}
\begin{equation} \label{eq:EigenVect0}
\Lambda_1 = \frac{1}{u_n} I,\Lambda_2=-\Lambda_1,
W_1  = A \Phi_1 + B \Phi_2=\frac{1}{\phi_0}\begin{pmatrix}  \phi_n \phi_0 \\ \phi_n^2 - (\tau_n - \phi_3) \\ \phi_n \phi_t - \tau_t \end{pmatrix},
W_2 = B \Phi_2=\frac{1}{\phi_0}\begin{pmatrix} 0 \\ - (\tau_n - \phi_3) \\ -\tau_t \end{pmatrix}.
\end{equation}
The details in (\ref{eq:enerBT5}) and (\ref{eq:EigenVect0}) reveal that five independent variables are involved in the boundary conditions.
\begin{remark}
Diagonalizing the boundary term (\ref{eq:enerBT3}) with standard eigenvalue techniques leads in general to very complex eigenvalues and eigenvectors, and hence complicated non-physical boundary conditions \cite{nordstrom2019}.
\end{remark}

\subsubsection{Boundary conditions of inflow-outflow type with nonzero external data}
We start by determining the number of boundary conditions \cite{nordstrom2020} required at inflow-outflow boundaries.
\begin{itemize}
\item At inflow $u_n < 0$, $W_1^T \Lambda_1 W_1 < 0$ and we need three conditions since $W_1$ has three components.
\item At outflow $u_n > 0$, $W_2^T \Lambda_2 W_2 < 0$ and we need two conditions since $W_2$ has two nonzero components.
\end{itemize}
Following \cite{nordstrom_roadmap}, the boundary conditions are applied weakly by inserting (\ref{eq:enerBT5}) and a lifting operator into (\ref{eq:ener3}):
\begin{align} \label{eq:ener4}
 \frac{d}{dt}  ||\Phi||_{\tilde I}^2 + & 2 \int_\Vol u_{i,j} \tau^*_{ij} d\Vol +  \oint\limits_{\partial\Omega}  \begin{pmatrix} W_1 \\ W_2 \end{pmatrix}^T \begin{pmatrix} \Lambda_1 & 0 \\  0 & \Lambda_2 \end{pmatrix} 
  \begin{pmatrix} W_1 \\ W_2 \end{pmatrix} ds + 2  \int_\Vol \Phi^T L(\Sigma \BC) d\Vol  = 0.
\end{align}
Here $\BC$ denotes the boundary conditions, $\Sigma$ is a penalty matrix and $L$ is a lifting operator implementing boundary conditions weakly. It is defined by $ \int_\Vol  \psi^T L(\phi) d\Vol = \oint\limits_{\partial\Omega} \psi^T \phi ds$,  where $\psi$ and $\phi$ are smooth vectors.

We will apply boundary conditions such that the boundary terms in (\ref{eq:ener4}) are bounded by external data only. 
Following \cite{NORDSTROM2024_BC}, we consider the 
nonlinear characteristic-like boundary condition on weak form
\begin{equation}\label{eq:BC-WEAk}
 \BC=\sqrt{|\Lambda^-|} W^- - {\bf{G}} = 0,
\end{equation}
where $\Lambda^-$ and $W^-$  are functions of the solution and $\bf{G}$ is external data.

To implement inflow boundary conditions weakly where $u_n < 0$, we need an operator $T_{1}$ so that
\begin{equation}
W^- = W_1 =  A \Phi_1 + B \Phi_2 =\frac{1}{\phi_0} \begin{pmatrix}  \phi_n \phi_0 \\ \phi_n^2 - (\tau_n - \phi_3) \\ \phi_n \phi_t - \tau_t \end{pmatrix} = T_1 \Phi = ( \Phi^T T_1^T)^T.
\end{equation}
\noindent At an outflow boundary where $u_n >0$ we require an operator $T_2$ so that
\begin{equation}
W^- =W_2 =  \\ B \Phi_2 =\frac{1}{\phi_0} \begin{pmatrix} 0 \\ - (\tau_n - \phi_3) \\ -\tau_t \end{pmatrix} = T_2 \Phi = ( \Phi^T T_2^T )^T.
\end{equation}
Implementing the boundary condition (\ref{eq:BC-WEAk}) weakly into (\ref{eq:ener4})  yield the augmented boundary term:
\begin{align}\label{eq:BT_SAT1}
\BT & = W^{+^T} \Lambda^+ W^- - W^{-^T} | \Lambda^- | W^- + 2\Phi^T  \Sigma \left( \sqrt{|\Lambda^-|}  W^- - {\bf{G}} \right). 
\end{align}

Next we choose $ \Sigma = T^T \sigma_0 \sqrt{|\Lambda^-|}$, where $T=T_1$ or $T_2$ depending on the signs of $ \Lambda_{1,2}$, which gives
\begin{align}\label{eq:BT_SAT2}
\BT & =  W^{+^T} \Lambda^+ W^- - (1- 2\sigma_0) W^{-^T} | \Lambda^- | W^- - 2 \sigma_0 W^{-^T} \sqrt{|\Lambda^-|} {\bf{G}}.
\end{align}
\noindent Choosing $\sigma_0 = 1$ as well as adding and subtracting ${\bf{G}}^T$ leads to an estimate in terms of data only since
\begin{align}\label{eq:BT_SAT3}
\BT  = 
 W^{+^T} \Lambda^+ W^- + \left( W^{-} \sqrt{|\Lambda^-|} - \bf{G} \right)^T \left( W^{-} \sqrt{|\Lambda^-|} - \bf{G} \right) - {\bf{G}}^T{\bf{G}}  \geq - {\bf{G}}^T{\bf{G}}.
\end{align}
\begin{remark}
The term related to $W^+$ in (\ref{eq:BT_SAT1}) is positive,  dissipative and requires no modification.
\end{remark}
\begin{remark}
IBVPs with estimates in terms of only external data are cited as strongly energy bounded  \cite{gustafsson1995time}.
\end{remark}

\subsubsection{Boundary conditions of solid wall  type with zero external data}
At solid wall boundaries, where $u_n=0$ we work directly on the $\BT$ in (\ref{eq:enerBT3}) which we group as follows
\begin{align} \label{eq:BCwall}
\BT = 
\left[ \Phi_1^T A \Phi_1 + 2 \Phi^T \Sigma_1 \BC_1 \right] + \left[ 2 \Phi_1^T B \Phi_2  + 2 \Phi^T \Sigma_2 \BC_2 \right].
\end{align}
We will remove the terms within brackets on the righthand-side of (\ref{eq:BCwall}) by appropriately selecting the penalty terms $\Sigma_1 \BC_1 $ and $\Sigma_2 \BC_2$.  Only the two boundary conditions $u_n = u_t = 0$ are available at a solid wall. For  consistency reasons they must cancel both bracketed terms in (\ref{eq:BCwall}) when imposed weakly.

For the first term in brackets in (\ref{eq:BCwall}), we set $\Sigma_1 = R_1^T \Sigma_1^\prime$ with $R_1 \Phi = \Phi_1$ to obtain
\begin{align}\label{SATwall1}
\Phi_1^T A \Phi_1 + 2 \Phi^T \Sigma_1 \BC_1 = \Phi_1^T A \Phi_1 + 2 \Phi_1^T \Sigma_1^\prime \BC_1,
\end{align}
\noindent where
\begin{equation}\label{Def_R1_N}
R_1 = \begin{pmatrix} 1 & 0 & 0 & 0 \\ 0 & n_1 & n_2 & 0 \\ 0 & -n_2 & n_1 & 0 \end{pmatrix}   \quad \text{and} \quad  \begin{pmatrix} n_1 & n_2 \\ -n_2 & n_1 \end{pmatrix}=N
\end{equation}
is the rotation matrix used in (\ref{eq:enerBT2}) above. Next, we insert  $\Sigma_1^\prime = \sigma_1$ and $\BC_1=A \Phi_1$ into (\ref{SATwall1}) to give
\begin{align}
\Phi_1^T A \Phi_1 + 2 \Phi_1^T \Sigma_1^\prime \BC_1 & = \Phi_1^T A \Phi_1 + 2 \Phi_1^T \left( \sigma_1 A \Phi_1 \right) = \Phi_1^T A \Phi_1 \left[1 + 2 \sigma_1 \right].
\end{align}
\noindent Equating $\sigma_1=-1/2$ yields stability. Consistency is proven by noting that $u_n$ is forced to zero by $\Sigma_1 \BC_1$,
since
\begin{equation*}
\Sigma_1 \BC_1 = \sigma_1 R_1^T A \Phi_1 = \sigma_1 R_1^T \Phi_1 u_n.
\end{equation*}
\noindent 

For the second term in brackets in (\ref{eq:BCwall}), we set $\Sigma_2 = R_2^T \Sigma_2^\prime$ with $R_2 \Phi = \Phi_2$ to obtain
\begin{align}\label{SATwall2}
2 \Phi_1^T B \Phi_2 + 2 \Phi^T \Sigma_2 \BC_2 &= 2 \left( B\Phi_2 \right)^T \Phi_1 + 2 \Phi^T_2 \Sigma_2^\prime \BC_2,
\end{align}
\noindent where $R_2= P_e T_r S_c$, $S_c = diag[1,1/\phi_0,1/\phi_0,1]$ and
\begin{equation}
P_e = \begin{pmatrix} 0 & 0 & 1 \\ 1 & 0 & 0 \\ 0 & 1 & 0 \end{pmatrix},~T_r = \begin{pmatrix} 0 & T_{r_{11}} & T_{r_{11}} & 0 \\ 0 & T_{r_{21}} & T_{r_{22}} & 0 \\ 0 & 0 & 0 & 1 \end{pmatrix}, ~ T_r = \mu N \left[ \begin{pmatrix} 2 n_1 & n_2 \\ 0 & n_1 \end{pmatrix} \pd{}{x} +  \begin{pmatrix} n_2 & 0 \\ n_1 & 2 n_2 \end{pmatrix} \pd{}{y}
\right].
\end{equation}
Setting $\Sigma_2^\prime=\sigma_2$ and $\BC_2=B^T \Phi_1$ into (\ref{SATwall2}) leads to,
\begin{align}\label{SATwall3}
2 \left( B\Phi_2 \right)^T \Phi_1 + 2 \Phi^T_2 \Sigma_2^\prime \BC_2 & = 2 \left( B\Phi_2 \right)^T \Phi_1 + 2 \sigma_2 \left( B \Phi_2\right)^T\Phi_1= 2 (1+\sigma_2) \left( B \Phi_2\right)^T\Phi_1.
\end{align}
Equating $\sigma_2=-1$ yields stability. Consistency  is proven by noting that $\Sigma_2 \BC_2$ forces $u_n$ and $u_t$ to zero, since
\begin{equation*}
\Sigma_2 \BC_2 = R_2^T \Sigma_2^\prime \BC_2  = \sigma_2 R_2^T B^T \Phi_1 =  \sigma_2 R_2^T B^T \begin{pmatrix} \phi_0  \\\phi_n \\ \phi_t \end{pmatrix} 
= \sigma_2 R_2^T \begin{pmatrix} u_n \\ -u_n \\ -u_t \end{pmatrix}.
\end{equation*}
\begin{remark} \label{Remtotal} 
The new inflow-outflow and solid wall boundary conditions leads to a strongly energy bounded VOF formulation with bounds on density, velocities and volume fraction in terms of external data only.
\end{remark}

\section{A stable straightforward nonlinear numerical approximation}\label{numerics}
The focus in this short note is on the continuous analysis, but for clarity we now sketch how the continuous formulation can be mimicked discretely, leading to 
a stable scheme.
The semi-discrete version of (\ref{eq:ge4}) (ignoring boundary conditions) using summation-by-parts (SBP) operators ${\bf D_{x_j}}$ \cite{nordstrom_roadmap,svard2014review,fernandez2014review} can be written
\begin{equation}\label{SWE_Disc}
{\bf \tilde I} \dfrac{d \vec \Phi}{dt}+ \frac{1}{2}  \Big[ {\bf D_{x_j}} ({\bf A_j}  \vec \Phi)+{\bf A_j} {\bf D_{x_j}}  \vec \Phi \Big] ={\bf B} {\bf D_{x_j}} {\bf S_j}.
\end{equation}
In (\ref{SWE_Disc}), ${\bf \tilde I}=  \tilde I \otimes I_x   \otimes I_y$ where $\otimes$ denotes the Kronecker product, the vector $\vec \Phi=(\vec \phi_0^T, \vec \phi_1^T,\vec \phi_2^T,\vec \phi_3^T)^T$ approximates $\Phi=(\phi_0, \phi_1,\phi_2,\phi_3)^T$ and the vector ${\bf S_j}=(0,\sec S_{1j}^T,\vec S_{2j}^T,\vec S_{3j}^T)^T$ approximates ${S_j}=(0,S_{1j} ,S_{2j},S_{3j})^T$ in each node.
The matrix elements of ${\bf A_j}$ and $ {\bf B} $ are matrices with node values of the matrix elements in $\tilde A_j$ and $\tilde  B$ injected on the diagonal as exemplified below in matrix ${\bf C}$
\begin{equation}
\label{illustration}
C =
\begin{pmatrix}
      c_{11}   &  \ldots  & c_{1n} \\
       \vdots   & \ddots & \vdots \\
       c_{n1} &  \ldots  & c_{nn}
\end{pmatrix}, \quad
{\bf C} =
\begin{pmatrix}
      {\bf c_{11}}   &  \ldots  &  {\bf c_{1n}}  \\
       \vdots          & \ddots &  \vdots           \\
        {\bf c_{n1}} &  \ldots &  {\bf c_{nn}} 
\end{pmatrix}, \quad
{\bf c_{ij}} =diag(c_{ij}(x_1,y_1), \ldots, c_{ij}(x_N,y_M)).
\end{equation}
Moreover, ${\bf D_{x_j}}=I_4  \otimes  \mathbb{D}_{x_j}$, $\mathbb{D}_{x_1}=D_x \otimes I_y$  and $\mathbb{D}_{x_2}=I_x \otimes D_y$ where 
$D_{x,y}$  are 1D SBP difference operators.
All matrices have appropriate sizes such that the matrix-matrix and matrix-vector operations are defined.

The discrete energy method (multiply (\ref{SWE_Disc}) from the left with  $2 \vec \Phi^T  {\bf P}$)  yields 
\begin{equation}\label{SWE_Disc_energy}
2  \vec \Phi^T {\bf P} \tilde {\bf I} \dfrac{d \vec \Phi}{dt}+ \vec \Phi^T ({\bf P} {\bf D_{x_j}} {\bf A_j} + {\bf A_j}{\bf P}{\bf D_{x_j}})\vec \Phi=2 \vec \Phi^T {\bf B} {\bf P}{\bf D_{x_j}} {\bf S_j},
\end{equation}
since ${\bf A_j}$ and ${\bf B}$ commute with  the diagonal symmetric positive definite integration operator ${\bf P}=I_4  \otimes \mathbb{P}$. Using the notation  $\| \vec \Phi \|_{\tilde {\bf P}}^2=\vec  \Phi^T {\bf P} \tilde {\bf I} \vec  \Phi$, noting that only the symmetric part of ${\bf P} {\bf D_{x_j}} {\bf A_j} + {\bf A_j}{\bf P}{\bf D_{x_j}}$ remains, and after applying the SBP properties (see  \cite{nordstrom_roadmap,svard2014review,fernandez2014review} for details) we obtain the semi-discrete energy rate
 \begin{equation}\label{Disc_energy_final}
 \dfrac{d}{dt} \| \vec  \Phi\|_{\tilde {\bf P}}^2 + 2 (\mathbb{D}_{x_j} \vec u_i)^T \mathbb{P} \vec{\mathbb{\tau^*}_{ij}} + \vec \Phi^T ({\bf B_{x_j}} {\bf A_j}) \vec \Phi-2 \vec u_i^T \mathbb{B}_{x_j} \vec{\mathbb{\tau}_{ij}} =0. 
\end{equation}
The semi-discrete energy rate (\ref{Disc_energy_final}) mimicks the continuous energy rate (\ref{eq:ener3}) perfectly on a rectangular domain. The second term from the left in (\ref{Disc_energy_final}) is the dissipation numerically integrated over the domain (using the volume integrator $\mathbb{P}$). The next two terms, contain the boundary terms numerically integrated along the boundary (using the boundary integrators $ {\bf B_{x_j}} =I_4 \otimes  \mathbb{B}_{x_j}$ and  $\mathbb{B}_{x_j}$).

\section{Summary and outlook}\label{sec:conclusion}
A new formulation of the incompressible multi-phase liguid-gas flow equations that leads to a well defined energy rate was derived. It was complemented with new weak boundary conditions of inflow-outflow and solid wall types that lead to a unique result for nonlinear VOF formulations: energy bounds of density, velocities and volume fraction in terms of external data only. 
The paper was concluded with a short illustration of how to construct a semi-discrete energy stable scheme by combining the new formulation with summation-by-parts operators. In future work we will develop nonlinear strongly energy stable schemes based on this new provably strongly energy bounded continuous formulation.


\section*{Acknowledgments}
\noindent J. N. was supported by Vetenskapsr{\aa}det, Sweden [award 2021-05484 VR] and University of Johannesburg Global Excellence and Stature Initiative Funding. A. G. M. was supported by the National Research Foundation (NRF) of South Africa  [grant 89916].  

\bibliographystyle{elsarticle-num}
\bibliography{References_Jan}

\end{document}